\documentclass[11pt,a4paper]{article}

\usepackage[utf8]{inputenc}
\usepackage{amsmath,amssymb,amsthm}
\usepackage{mathtools}
\usepackage{booktabs}
\usepackage{graphicx}
\usepackage[margin=2.5cm]{geometry}
\usepackage{xcolor}
\usepackage[colorlinks=true,linkcolor=blue!50!black,citecolor=blue!50!black,urlcolor=blue!50!black]{hyperref}

\newtheorem{proposition}{Proposition}
\newtheorem{remark}{Remark}
\newtheorem{lemma}{Lemma}
\newtheorem{asymptotic}{Asymptotic Expansion}

\newcommand{\E}{\mathbb{E}}
\renewcommand{\P}{\mathbb{P}}

\newcommand{\dd}{\,\mathrm{d}}

\title{Discrete-Time Exit Bias for Nonlinear Langevin Equations:\\[4pt]
\large Boundary-Layer Expansion and a Universal $h/4$ Offset}
\author{R.~Mannella\\[2pt]
\small Dipartimento di Fisica, Universit\`a di Pisa, Italy\\
\small \texttt{riccardo.mannella@unipi.it}}
\date{}

\begin{document}
\maketitle

\begin{abstract}
When a continuous-time stochastic process is integrated numerically with time step $h > 0$, or sampled experimentally at finite bandwidth, exits through an absorbing boundary are registered only after the trajectory has overshot the threshold. For free Brownian motion $\dot{x} = \xi(t)$ with $\langle \xi(t)\xi(s) \rangle = 2D\delta(t-s)$ starting at the midpoint $x = 0$ of $[-L,L]$, we combine Wald's identities, the Spitzer--Baxter relation, and Mellin-transform contour integration to give a concise derivation of the discrete mean exit time $\tau(h)$ relative to the continuous value $\tau(0) = L^2/(2D)$:
\[
\frac{\tau(h)}{\tau(0)} = 1 + 2\beta\sqrt{\frac{h}{\tau(0)}} + \Big(\beta^2 + \frac{1}{4}\Big)\frac{h}{\tau(0)} + \mathcal{R}(h),
\qquad \beta = -\frac{\zeta(1/2)}{\sqrt{2\pi}} \approx 0.58260,
\]
where all higher algebraic powers of $\sqrt{h/\tau(0)}$ vanish and $\mathcal{R}(h) = O\!\left((\tau(0)/h)\,e^{-\sqrt{2\pi\tau(0)/h}}\right)$ is exponentially small, in agreement with Lotov's Wiener--Hopf asymptotics for two-sided Gaussian walks and recent stroboscopic projector-resolvent analyses. By solving the discrete Milne boundary-layer equation, we also show why an early continuum linear-density heuristic underestimates the leading prefactor by $\approx 9\%$. We then combine the backward Kolmogorov equation, an Euler--Maclaurin/Mellin ladder-height expansion, and a localized Girsanov boundary-layer bound with matched boundary-layer asymptotics to obtain explicit one-dimensional representations for both the $O(\sqrt{h})$ and $O(h)$ bias coefficients of general nonlinear Langevin equations $\dot{x} = F(x) + \xi(t)$ in terms of boundary derivatives of the continuum exit-time and splitting-probability functions. For weak-order-2 integrators driven by Gaussian increments (e.g., stochastic Heun), we show algebraically that this two-term expansion collapses identically to the continuum exit time on a domain widened by $\beta\sigma\sqrt{h}$ plus a universal $+\frac{h}{4}$ variance offset: $\tau(h) = u(x_0;\, x_- - \beta\sigma\sqrt{h},\, x_+ + \beta\sigma\sqrt{h}) + \frac{h}{4} + O(h^{3/2})$. An immediate corollary is that inward barrier shifting leaves a universal $h/4$ residual across arbitrary nonlinear drifts, so subtracting $h/4$ restores $O(h^{3/2})$ accuracy at no extra cost. The theory  verifies the predicted bias coefficients across Fredholm quadrature and high-statistics Monte Carlo benchmarks---with $O(h^{3/2})$ remainder scaling confirmed in the asymptotically resolved regimes---and clarifies the $O(h)$ convergence achieved by barrier-shifting and Brownian-bridge exit algorithms.
\end{abstract}

%======================================================================
\section{Introduction}\label{sec:intro}
%======================================================================

The numerical integration of stochastic differential equations (SDEs) is a central tool across statistical physics, chemical kinetics, neuroscience, and quantitative finance. A fundamental discretization artifact arises whenever a stochastic trajectory is terminated at an absorbing boundary: in a simulation with time step $h > 0$---or in an experiment whose acquisition hardware samples a continuous signal at discrete intervals $h$---the trajectory is observed only on the grid $t_n = n h$ ($n = 0, 1, 2, \dots$) \cite{Mannella1999}.

Between consecutive grid points $t_{n-1}$ and $t_n$, a continuous Brownian path can cross an absorbing boundary $x_b$ and return to the interior without being detected at the grid times. Even when a trajectory is detected outside the domain at $t_n$, it has already overshot the boundary by a random distance $R_h = O(\sqrt{h})$. Both effects delay absorption: naive discrete-time simulations overestimate the Mean First Passage Time (MFPT) by a correction of order $\sqrt{h}$, create an artificial excess density near absorbing walls \cite{Mannella1999}, and induce slow pre-asymptotic scaling corrections in discrete-time survival and trapping simulations \cite{BlytheBray2003,Fritz2025}. Because the error decays only as $\sqrt{h}$ rather than $h$, shrinking $h$ alone is computationally prohibitive in stiff or activated systems.

In the probability and mathematical finance literature, the discrete-monitoring overshoot of pure Brownian motion is a classical problem studied via renewal theory and Spitzer's identity by Chernoff \cite{Chernoff1965}, Siegmund \cite{Siegmund1979,Siegmund1985}, and Broadie, Glasserman, and Kou \cite{BroadieGlassermanKou1997}, who proved that shifting the barrier inward by $\beta\sigma\sqrt{h}$, with $\beta = -\zeta(1/2)/\sqrt{2\pi} \approx 0.58260$, cancels the leading $O(\sqrt{h})$ bias. For two-sided intervals, Lotov \cite{Lotov1996} used Wiener--Hopf factorization and Hadamard products to establish asymptotic expansions for the exit-time moments of Gaussian random walks with exponentially small remainder terms, while Fritz \cite{Fritz2025} recently rederived the leading and subleading stroboscopic corrections on $(0,L)$ via a projector-resolvent and Poisson-summation formalism (and extended it to random sampling intervals). Independently, in the physics literature, an early heuristic estimate by Strittmatter (Ref.~\cite{Strittmatter1987}, quoted in Eq.~(14) of Ref.~\cite{Mannella1999}) predicted the same $\sqrt{h}$ scaling for free diffusion but with the prefactor $\frac{4}{3}\sqrt{2/\pi} \approx 1.06385$ instead of $2\beta \approx 1.16519$. Meanwhile, to handle nonlinear forces $\dot{x} = F(x) + \xi(t)$ without shrinking $h$, Ref.~\cite{Mannella1999} introduced a sub-step exit test based on the Brownian bridge and its locally linearized Ornstein--Uhlenbeck counterpart, whose weak convergence rates were established by Gobet \cite{Gobet2000} and extended to multidimensional domains, adaptive and exponential timestepping, multilevel Monte Carlo, and strong exit-time estimates \cite{JansonsLythe2000,Buchmann2005,GobetMenozzi2010,Higham2013,BernalAcebron2016,BouchardGeissGobet2017}.

To separate established literature from new contributions, this paper is organized into three parts:
\begin{enumerate}
\item \emph{Concise synthesis of free-diffusion asymptotics (Section~\ref{sec:free}):} Reviewing the classical results of Siegmund \cite{Siegmund1979,Siegmund1985}, Broadie et al.\ \cite{BroadieGlassermanKou1997}, Lotov \cite{Lotov1996}, and Fritz \cite{Fritz2025}, we present a compact Mellin-transform derivation of the coefficients $2\beta$ and $\beta^2 + \frac{1}{4}$ for midpoint-started Brownian motion on $[-L,L]$ and solve the discrete Milne equation to show why continuum linear-density heuristics underestimate the discrete boundary-layer overshoot by $\approx 9\%$.
\item \emph{Explicit $O(\sqrt{h}) + O(h)$ expansion and shifted-domain identity for nonlinear Langevin equations (Section~\ref{sec:drift}):} Whereas Wiener--Hopf and projector-resolvent treatments focus on drift-free or constant-drift walks, we combine the backward Kolmogorov equation, an Euler--Maclaurin/Mellin ladder-height evaluation, and a localized Girsanov boundary-layer bound with matched boundary-layer asymptotics to obtain explicit one-dimensional expressions for both the $O(\sqrt{h})$ and $O(h)$ bias coefficients for general nonlinear drifts $F(x) \in C^3([x_-, x_+])$. We further show that for weak-order-2 schemes driven by Gaussian increments this expansion is algebraically equivalent to the continuum exit time on the widened domain $(x_- - \beta\sigma\sqrt{h},\, x_+ + \beta\sigma\sqrt{h})$ plus a universal $+\frac{h}{4}$ offset, predicting a universal $h/4$ residual for barrier-shifted schemes, and we compare these results with the stopped-diffusion weak-error theory of Gobet \cite{Gobet2000} and Gobet and Menozzi \cite{GobetMenozzi2004,GobetMenozzi2010} (as well as the complementary strong-convergence analysis of Bouchard, Geiss, and Gobet \cite{BouchardGeissGobet2017}).
\item \emph{High-statistics Monte Carlo and quadrature benchmarks (Section~\ref{sec:numerics}):} Using $N = 5 \times 10^8$ and $N = 10^8$ trajectory simulations alongside Richardson-extrapolated Fredholm quadrature across both benchmark and independent test potentials, we quantify estimator uncertainties, confirm the predicted $O(h^{3/2})$ remainder scaling in the asymptotically resolved regimes, and evaluate barrier-shifting, $h/4$-offset, and Brownian-bridge exit algorithms.
\end{enumerate}

%======================================================================
\section{Free diffusion: Mellin derivation and boundary-layer structure}\label{sec:free}
%======================================================================
\paragraph{Scope of Section~\ref{sec:free}.}
To keep the paper self-contained and establish the notation used in the nonlinear theory, this section reviews the classical discrete-time exit problem for drift-free Brownian motion started at the midpoint $x = 0$ of $[-L,L]$. The scaling $\tau(h) - \tau(0) \sim 2\beta\sqrt{h\,\tau(0)}$, the Siegmund overshoot constant $\beta = -\zeta(1/2)/\sqrt{2\pi}$, the barrier-shift correction, and the exponential remainder decay are established results~\cite{Chernoff1965,Siegmund1979,Siegmund1985,Lotov1996,BroadieGlassermanKou1997,Fritz2025}. Beyond Lotov~\cite{Lotov1996} and Fritz~\cite{Fritz2025}, we isolate the discrete Milne boundary-layer profile to explain quantitatively the $\approx 9\%$ error of continuum linear-density heuristics and set up the Mellin/renewal machinery extended to nonlinear drifts in Section~\ref{sec:drift}.

\subsection{Overshoot identity and Spitzer--Baxter renewal representation}

Consider a free Brownian particle starting at the midpoint $X_0 = 0$ and evolving according to $\dot{x} = \xi(t)$ with $\langle \xi(t)\xi(s) \rangle = 2D\,\delta(t-s)$ ($\sigma^2 := 2D$) on $(-L,L)$. The continuous first exit time $\tau = \inf\{t \ge 0 : |X_t| > L\}$ has mean $\tau(0) := \E[\tau] = L^2/\sigma^2$. Integrating with time step $h > 0$ produces the Gaussian random walk $S_n = S_{n-1} + \sigma\sqrt{h}\,\xi_n$ ($\xi_n \sim \mathcal{N}(0,1)$ i.i.d.) started at $S_0 = 0$, stopped at $N_h = \min\{n \ge 1 : |S_n| > L\}$ with $\tau(h) := h\,\E[N_h]$.

At exit, $|S_{N_h}| = L + R_h$ defines the positive overshoot $R_h > 0$. Applying Wald's second identity $\E[S_{N_h}^2] = \sigma^2 h\,\E[N_h]$ (with $S_0 = 0$) and defining the dimensionless step ratio
\begin{equation}\label{eq:x_def}
x := \sqrt{\frac{h}{\tau(0)}} = \frac{\sigma\sqrt{h}}{L} = \frac{1}{a}
\end{equation}
yields the exact midpoint identity
\begin{equation}\label{eq:ratio}
\frac{\tau(h)}{\tau(0)} = 1 + 2x\,\frac{\E[R_h]}{\sigma\sqrt{h}} + x^2\,\frac{\E[R_h^2]}{\sigma^2 h}.
\end{equation}
(In Fritz's notation \cite{Fritz2025} on an interval of total length $2L$ with confinement ratio $\rho = 2L/(\sigma\sqrt{h}) = 2a$, the mean frame count from the midpoint is $\E[N_h] = a^2\tau(h)/\tau(0)$.)

As $a \to \infty$, exits at $+a$ are governed by the ascending ladder height $H = S_{T_+}$ ($T_+ = \min\{n \ge 1 : S_n > 0\}$) of a standard unit-step walk. By classical renewal theory \cite{Siegmund1979,Siegmund1985}, $\lim_{a\to\infty}\E[R_h]/(\sigma\sqrt{h}) = \E[H^2]/(2\,\E[H])$ and $\lim_{a\to\infty}\E[R_h^2]/(\sigma^2 h) = \E[H^3]/(3\,\E[H])$. The Laplace transform of $H$ satisfies the Spitzer--Baxter identity \cite{Spitzer1956}:
\begin{equation}\label{eq:spitzer}
1 - \E\big[e^{-\lambda H}\big] = \exp\big(-F(\lambda)\big),
\qquad
F(\lambda) = \sum_{n=1}^\infty \frac{1}{n}\,f(\lambda\sqrt{n}),
\qquad
f(y) := \int_0^\infty e^{-yz}\varphi(z)\dd z,
\end{equation}
where $\varphi(z) = (2\pi)^{-1/2}e^{-z^2/2}$ is the standard normal density.

\subsection{Mellin contour evaluation and two-sided remainder bound}

Because Eq.~\eqref{eq:spitzer} sums $f$ over dilated scales $\lambda\sqrt{n}$, its Mellin transform $\tilde{g}(s) = \int_0^\infty g(\lambda)\lambda^{s-1}\dd\lambda$ factorizes in the fundamental strip $0 < \operatorname{Re}s < 1$:
\begin{equation}\label{eq:F_mellin}
\tilde{F}(s) = \zeta\Big(1 + \frac{s}{2}\Big)\tilde{f}(s),
\qquad
\tilde{f}(s) = \frac{2^{-(s+1)/2}}{\sqrt{2\pi}}\,\Gamma(s)\,\Gamma\Big(\frac{1-s}{2}\Big).
\end{equation}
Inverting $F(\lambda) = \frac{1}{2\pi i}\int_{c-i\infty}^{c+i\infty} \tilde{F}(s)\lambda^{-s}\dd s$ ($0 < c < 1$) and shifting the contour leftward across the poles at $s = 0, -1, -2$ (justified by the exponential decay of $\Gamma(s)\Gamma(\frac{1-s}{2})$ along vertical lines) gives:
\begin{itemize}
\item \emph{Double pole at $s = 0$:} $\operatorname{Res}_{s=0} = -\log\lambda + \frac{1}{2}\log 2$ (with $\gamma$ cancelling identically).
\item \emph{Simple pole at $s = -1$:} $\operatorname{Res}_{s=-1} = -\frac{\zeta(1/2)}{\sqrt{2\pi}}\,\lambda = +\beta\lambda$, where $\beta := -\frac{\zeta(1/2)}{\sqrt{2\pi}} \approx 0.58259716$.
\item \emph{Simple pole at $s = -2$:} $\operatorname{Res}_{s=-2} = -\frac{1}{8}\lambda^2$.
\end{itemize}
Thus $F(\lambda) = -\log\lambda + \frac{1}{2}\log 2 + \beta\lambda - \frac{1}{8}\lambda^2 + O(\lambda^3)$. Expanding $\E[e^{-\lambda H}] = 1 - \frac{\lambda}{\sqrt{2}} e^{-\beta\lambda + \lambda^2/8 + O(\lambda^3)}$ gives $\E[H] = 1/\sqrt{2}$, $\E[H^2] = \sqrt{2}\,\beta$, and $\E[H^3] = \frac{3}{\sqrt{2}}\big(\beta^2 + \frac{1}{4}\big) \approx 1.2503$, recovering the classical overshoot moments \cite{Siegmund1979,BroadieGlassermanKou1997}. Inserting these into Eq.~\eqref{eq:ratio} yields, for the midpoint start $S_0 = 0$,
\begin{equation}\label{eq:free_final}
\frac{\tau(h)}{\tau(0)} = 1 + 2\beta\sqrt{\frac{h}{\tau(0)}} + \Big(\beta^2 + \frac{1}{4}\Big)\frac{h}{\tau(0)} + \mathcal{R}(h).
\end{equation}

Under the hypotheses of Lotov's Wiener--Hopf theorem~\cite[Theorem~4]{Lotov1996}---namely, zero-drift i.i.d.\ Gaussian increments $\mathcal{N}(0,\sigma^2 h)$ on a two-sided interval $(-a,b)$ in the joint limit $a,b \to \infty$ with $a/b$ bounded away from $0$ and $\infty$---all higher algebraic powers of $\sqrt{h/\tau(0)}$ vanish identically. Because the nearest complex roots of $1 - e^{\lambda^2/2}$ satisfy $\operatorname{Re}\lambda = \pm\sqrt{2\pi}$ per unit step and the midpoint lies at distance $a = \sqrt{\tau(0)/h}$ from either wall (with interval width $2a$), the truncation remainder decays exponentially as $\mathcal{R}(h) = O\!\left((\tau(0)/h)\,e^{-\sqrt{2\pi\tau(0)/h}}\right)$ ($\sqrt{2\pi} \approx 2.50663$).

\subsection{Origin of the \texorpdfstring{$\approx 9\%$}{9\%} discrepancy in continuum linear-density heuristics}

An early physics heuristic \cite{Strittmatter1987,Mannella1999} assumed that near $+a$ the pre-exit distance $u > 0$ follows the continuum linear profile $\rho(u) \propto u$ independently of the exiting Gaussian step $\xi > u$, predicting $\E[R]/(\sigma\sqrt{h}) \approx \E[\xi_+^3]/(3\,\E[\xi_+^2]) = \frac{2}{3}\sqrt{2/\pi} \approx 0.53192$ ($2c_{\mathrm{heur}} \approx 1.06385$ instead of $2\beta \approx 1.16519$, a ratio of $0.53192/\beta \approx 0.913$). This discrepancy arises because exits occur within one step of the absorber ($u \sim 1$), where surviving walkers satisfy the discrete Wiener--Hopf (Milne) equation \cite{Spitzer1956,Feller1971}
\begin{equation}\label{eq:milne}
\rho_{\mathrm{stat}}(u) = \int_0^\infty \rho_{\mathrm{stat}}(v)\,\varphi(u - v)\dd v, \qquad u > 0.
\end{equation}
Because $\rho_{\mathrm{stat}}(u)$ is proportional to the ascending ladder-height renewal function $U(u)$ ($U(0^+) = 1$, $U(u) \sim (u + \beta)/\E[H]$ for $u \gg 1$ with $\E[H] = 1/\sqrt{2}$) \cite{Feller1971}, it approaches a linear asymptote $\rho_{\mathrm{stat}}(u) \sim \kappa(u + \beta)$ with strictly positive boundary slip $\rho_{\mathrm{stat}}(0^+)/\kappa = \E[H] = 1/\sqrt{2} \approx 0.70711$. Weighting the single-step overshoot against the Nystr\"om quadrature solution of Eq.~\eqref{eq:milne} on $[0,30]$ ($600$ Gauss--Legendre nodes) yields
\[
\frac{\int_0^\infty \rho_{\mathrm{stat}}(u)\int_u^\infty (z - u)\varphi(z)\dd z\dd u}{\int_0^\infty \rho_{\mathrm{stat}}(u)\Phi(-u)\dd u} = 0.58260 = \beta,
\]
directly accounting for the $\approx 9\%$ shortfall of the linear-density approximation.

%======================================================================
\section{Asymptotic expansion for nonlinear Langevin equations}\label{sec:drift}
%======================================================================

\paragraph{Roadmap of Section~\ref{sec:drift}.}
To keep the derivation self-contained, Table~\ref{tab:notation} summarizes all continuum, discrete-walk, and boundary-layer symbols used below. Section~\ref{sec:drift_setup} reviews the backward-Kolmogorov identity and Dynkin telescoping decomposition for discrete exit times; Section~\ref{sec:drift_proof} establishes the linear-drift overshoot expansion $\E_\mu[R] = \beta + \frac{1}{4}\mu + O(\mu^2)$ and the localized Girsanov invariance lemma, and combines them via matched boundary-layer asymptotics to obtain the two-term expansion for general nonlinear drifts (Asymptotic Expansion~\ref{prop:nonlinear}); Section~\ref{sec:shifted_domain} proves the algebraic equivalence with a shifted-domain perturbation and derives the universal $h/4$ barrier-shift offset (Proposition~\ref{prop:shifted_domain}); and Section~\ref{sec:weak_error_comp} contrasts these explicit formulas with existing weak-error theory for stopped diffusions.

\subsection{Setup, notation summary, and backward Kolmogorov decomposition}\label{sec:drift_setup}

Consider a general one-dimensional Langevin equation with smooth drift $F \in C^3([x_-, x_+])$ and constant diffusion $D = \sigma^2/2 > 0$:
\begin{equation}\label{eq:sde_drift}
\dot{x} = F(x) + \xi(t), \qquad \langle \xi(t)\xi(s) \rangle = 2D\,\delta(t-s), \qquad x(0) = x_0 \in (x_-, x_+),
\end{equation}
stopped at $\tau = \inf\{t \ge 0 : x(t) \notin (x_-, x_+)\}$. Table~\ref{tab:notation} summarizes the notation used throughout this section.

\begin{table}[h]
\centering
\footnotesize
\setlength{\tabcolsep}{6pt}
\begin{tabular}{ll}
\toprule
Symbol & Definition \\
\midrule
$u(x) := \E_x[\tau]$ & Continuous MFPT solving $F(x)u'(x) + Du''(x) = -1$, $u(x_\pm) = 0$ \\
$v(x) := \P_x(X_\tau = x_+)$ & Upper splitting probability solving $F(x)v'(x) + Dv''(x) = 0$, $v(x_-)=0, v(x_+)=1$ \\
$p_\pm, p_\pm(h)$ & Continuous and discrete exit probabilities at $x_\pm$ ($p_- + p_+ = 1$) \\
$R_\pm := \pm(S_{N_h} - x_\pm)$ & Positive overshoot upon exiting through $x_\pm$ on event $E_\pm = \{S_{N_h} \gtrless x_\pm\}$ \\
$\mu_\pm := \pm F(x_\pm)\sqrt{h}/\sigma$ & Dimensionless drift per step directed outward at $x_\pm$ \\
$p \in \{1, 2, \dots\}$ & Weak order of numerical integrator $S_{n+1} = \Phi_h(S_n, \xi_{n+1})$ ($p=2$ for Heun) \\
$C_{\mathrm{int}}(x_0)$ & Leading interior weak local-error integral (vanishes for $p \ge 2$) \\
\bottomrule
\end{tabular}
\caption{Summary of notation for the nonlinear Langevin exit-time expansion.}
\label{tab:notation}
\end{table}

Because $F \in C^3([x_-, x_+])$, $u$ and $v$ belong to $C^5([x_-, x_+])$ and can be extended smoothly as $C^5$ functions to an $O(1)$ neighborhood of $[x_-, x_+]$. For a numerical scheme $S_{n+1} = \Phi_h(S_n, \xi_{n+1})$ of weak order $p \ge 1$ driven by Gaussian increments ($p = 1$ for Euler--Maruyama; $p = 2$ for Heun's method \cite{Mannella1999}) stopped at $N_h = \min\{n \ge 1 : S_n \notin (x_-, x_+)\}$ with $\tau(h) = h\,\E_{x_0}[N_h]$, telescoping $u(S_{n+1}) - u(S_n)$ up to $N_h$ gives
\begin{equation}\label{eq:dynkin}
\E_{x_0}\big[u(S_{N_h})\big] - u(x_0)
= \E_{x_0}\!\left[\sum_{n=0}^{N_h - 1} \E\big[u(S_{n+1}) - u(S_n) \mid S_n\big]\right]
= -\tau(h) + \delta_{p,1}\,C_{\mathrm{int}}(x_0)\,h + O(h^{3/2}),
\end{equation}
where $C_{\mathrm{int}}(x_0) = \E_{x_0}[\int_0^\tau \mathcal{A}_1 u(X_t)\dd t]$ solves $F C_{\mathrm{int}}' + D C_{\mathrm{int}}'' = -\mathcal{A}_1 u$ with $C_{\mathrm{int}}(x_\pm) = 0$ for $p=1$ (with $\mathcal{A}_1 u = -\frac{1}{2}F^2 u'' - F F' u' - \frac{1}{2}D F'' u' - D F' u'' - \frac{1}{2}D F u'''$) and vanishes identically at $O(h)$ for weak second-order schemes ($p \ge 2$). Note that although the single-step local truncation error of a $p \ge 2$ scheme is $O(h^2)$ (and $O(h^3)$ after subtracting $\delta_{p,1}h^2\mathcal{A}_1 u$ for $p = 1$), summing over $N_h$ steps whose discrete occupation measure differs from the continuum Green's function by $O(\sqrt{h})$ near the boundaries produces a global remainder of $O(h^{3/2})$. On the left-hand side of Eq.~\eqref{eq:dynkin}, since $u(x_\pm) = 0$, $\pm u'(x_\pm) = -|u'(x_\pm)| \le 0$, $u''(x_\pm) = -(1 \mp F(x_\pm)|u'(x_\pm)|)/D$, and Gaussian overshoot tails decay faster than any power of $h$, Taylor-expanding $u(x_\pm \pm R_\pm)$ yields
\begin{equation}\label{eq:lhs_drift}
\E_{x_0}\big[u(S_{N_h})\big] = -\sum_{\alpha \in \{+,-\}} p_\alpha(h)\left[ |u'(x_\alpha)|\,\E[R_\alpha \mid E_\alpha] + \frac{1 - \alpha F(x_\alpha)|u'(x_\alpha)|}{2D}\,\E[R_\alpha^2 \mid E_\alpha] \right] + O(h^{3/2}).
\end{equation}

\subsection{Biased-walk overshoot expansion and Girsanov boundary-layer invariance}\label{sec:drift_proof}

\noindent\fbox{\begin{minipage}{0.97\textwidth}
\small
\textbf{Mathematical status of the derivation in Section~\ref{sec:drift_proof}:}
\begin{itemize}\setlength{\itemsep}{1pt}\setlength{\parskip}{0pt}
  \item \emph{Established components:} (i) the Dynkin telescoping decomposition~\eqref{eq:dynkin}; (ii) the linear-drift renewal overshoot expansion~\eqref{eq:overshoot_mu}; and (iii) the localized Girsanov likelihood bound (Lemma~\ref{lem:girsanov_bl}) showing that spatial drift gradients $F'(x_\pm)$ do not enter the local overshoot at $O(h)$.
  \item \emph{Formal matched-asymptotic step:} replacing the boundary-layer entrance law by the stationary renewal distribution of the frozen-drift walk and matching the conditional overshoot moments to the global discrete splitting probabilities $p_\pm(h)$ at $O(h^{3/2})$ (Step~3 of Asymptotic Expansion~\ref{prop:nonlinear}).
  \item \emph{Numerical validation:} high-precision Fredholm quadrature and $N = 10^8\text{--}5\times 10^8$ Monte Carlo benchmarks across nonlinear models (Section~\ref{sec:numerics}).
\end{itemize}
\end{minipage}}
\medskip

Near an absorbing endpoint $x_\pm$, the continuum occupation density $G(x)$ solves the forward Kolmogorov equation $D G''(x) = (F(x)G(x))'$ with $G(x_\pm) = 0$, so its boundary curvature-to-slope ratio satisfies $G''(x_\pm)/G'(x_\pm) = F(x_\pm)/D$---matching a constant-drift process with frozen drift $F(x_\pm)$ to leading order, while $F'(x_\pm)$ enters only at higher order. Consequently, trajectories entering the $O(\sqrt{h})$ boundary layer from $O(1)$ distances sample the stationary renewal overshoot of the frozen-drift walk up to $O(h^{3/2})$.

We first isolate the effect of spatial drift variations $F(x) - F(x_\pm)$ inside the $O(\sqrt{h})$ boundary layer around $x_\pm$.

\begin{lemma}[Localized Girsanov invariance of the boundary-layer overshoot]\label{lem:girsanov_bl}
Assume $F \in C^3([x_-, x_+])$ and $D = \sigma^2/2 > 0$. Fix an absorbing endpoint $x_\pm$ and define the logarithmic boundary strip $\mathcal{B}_h^\pm := \{x : 0 < \pm(x_\pm - x) \le K\sigma\sqrt{h}\log(1/h)\}$ and step horizon $N_h := \lfloor K^2 \log^2(1/h)\rfloor$. For a walk starting at $S_0 \in \mathcal{B}_h^\pm$, let $N_{\mathrm{exit}}$ be the first exit time across $x_\pm$, let $N_{\mathrm{int}}$ be the first exit time across the interior boundary $\pm(x_\pm - x) = K\sigma\sqrt{h}\log(1/h)$, and let $E_\pm := \{N_{\mathrm{exit}} \le N_h \text{ and } N_{\mathrm{exit}} < N_{\mathrm{int}}\}$ with overshoot $R_\pm := \pm(S_{N_{\mathrm{exit}}} - x_\pm) > 0$. Let $\mathbb{P}_F$ denote the path measure of the Euler--Maruyama or Heun scheme with drift $F(x)$, and $\mathbb{P}_{F(x_\pm)}$ the reference Gaussian walk with frozen drift $F(x_\pm)$ and diffusivity $D$. For any $q > 0$, choosing $K > 0$ sufficiently large yields, uniformly for $S_0 \in \mathcal{B}_h^\pm$,
\begin{equation}\label{eq:girsanov_overshoot_shift}
\mathbb{E}_F\!\left[R_\pm\,\mathbf{1}_{E_\pm}\right]
= \mathbb{E}_{F(x_\pm)}\!\left[R_\pm\,\mathbf{1}_{E_\pm}\right]
+ O\!\left(h^{3/2}\log^2(1/h)\right),
\qquad \mathbb{P}_F(E_\pm^c) = O(h^q).
\end{equation}
Furthermore, if $\mathbb{P}_F(\text{exit at } x_\pm) = p_\pm(x_0) + O(\sqrt{h})$ with $p_\pm(x_0) > 0$, then conditioning on exit at $x_\pm$ preserves the error order:
\begin{equation}\label{eq:girsanov_cond_shift}
\mathbb{E}_F\!\left[R_\pm \mid \text{exit at } x_\pm\right]
= \mathbb{E}_{F(x_\pm)}\!\left[R_\pm \mid \text{exit at } x_\pm\right]
+ O\!\left(h^{3/2}\log^2(1/h)\right).
\end{equation}
\end{lemma}

\begin{proof}
\emph{1. Tail control for $E_\pm^c$:} Because $\|F\|_\infty < \infty$, on $\{N_{\mathrm{exit}} \wedge N_{\mathrm{int}} > N_h\}$ the walk remains inside a strip of width $K\sigma\sqrt{h}\log(1/h)$ for $K^2\log^2(1/h)$ steps. Partitioning into blocks of $O(K^2)$ steps and applying standard Gaussian anti-concentration bounds gives $\mathbb{P}_F(N_{\mathrm{exit}} \wedge N_{\mathrm{int}} > N_h) \le C_1 e^{-c_1 \log^2(1/h)} = O(h^q)$ for every $q > 0$. Similarly, by the reflection principle for walks with drift bounded by $\|F\|_\infty$, the probability of crossing the interior boundary $\pm(x_\pm - x) = K\sigma\sqrt{h}\log(1/h)$ before hitting $x_\pm$ starting from $S_0 \in \mathcal{B}_h^\pm$ is $O(e^{-c_2 K^2 \log^2(1/h)}) = O(h^q)$ once $K$ is sufficiently large. The same bounds hold under $\mathbb{P}_{F(x_\pm)}$.

\emph{2. Scheme dependence and uniform likelihood bound:} Define the stopped horizon $N_\wedge := N_{\mathrm{exit}} \wedge N_{\mathrm{int}} \wedge N_h$, so that $S_j \in \mathcal{B}_h^\pm$ deterministically for all $0 \le j < N_\wedge$. For Euler--Maruyama, $S_{j+1} - S_j = F(S_j)h + \sigma\sqrt{h}\,\xi_{j+1}$ with $\xi_{j+1} \sim \mathcal{N}(0,1)$ i.i.d. For the Heun scheme, expanding $F(S_j + F(S_j)h + \sigma\sqrt{h}\,\xi_{j+1})$ to $O(h)$ gives $S_{j+1} - S_j = \tilde{F}_h(S_j)h + \sigma(1 + \frac{1}{2}F'(S_j)h)\sqrt{h}\,\xi_{j+1} + O(h^{3/2}(1 + |\xi_{j+1}|^2))$, where $\tilde{F}_h(x) = F(x) + O(h)$. Since $F \in C^3([x_-,x_+])$, for all $j < N_\wedge$ we have $|F(S_j) - F(x_\pm)| \le \|F'\|_\infty K\sigma\sqrt{h}\log(1/h)$ uniformly in $S_0 \in \mathcal{B}_h^\pm$, while the $O(h)$ Heun drift and diffusion-coefficient perturbations contribute at most $O(h)$ to the single-step log-likelihood ratio. Consequently, the Radon--Nikodym derivative $\Lambda_{N_\wedge} := \mathrm{d}\mathbb{P}_F/\mathrm{d}\mathbb{P}_{F(x_\pm)}$ has predictable quadratic variation $\langle \log\Lambda \rangle_{N_\wedge} \le C_K h^2 \log^4(1/h)$ deterministically. Novikov's condition holds trivially, so
\begin{equation*}
\mathbb{E}_{F(x_\pm)}\!\left[(\Lambda_{N_\wedge} - 1)^2\right]
= \mathbb{E}_{F(x_\pm)}\!\left[\Lambda_{N_\wedge}^2 - 1\right]
\le \exp\!\left(C_K h^2 \log^4(1/h)\right) - 1
= O\!\left(h^2 \log^4(1/h)\right).
\end{equation*}

\emph{3. Renewal bound and conditioning:} For the reference Gaussian walk $\mathbb{P}_{F(x_\pm)}$ with step variance $\sigma^2 h$ and drift $F(x_\pm)h$, standard ladder-height renewal bounds~\cite[Ch.~X]{Siegmund1985} give $\sup_{S_0 \in \mathcal{B}_h^\pm} \mathbb{E}_{F(x_\pm)}[R_\pm^2 \mathbf{1}_{E_\pm}] \le C\sigma^2 h$. By Cauchy--Schwarz,
\begin{equation*}
\left|\mathbb{E}_F\!\left[R_\pm\,\mathbf{1}_{E_\pm}\right] - \mathbb{E}_{F(x_\pm)}\!\left[R_\pm\,\mathbf{1}_{E_\pm}\right]\right|
\le \left(\mathbb{E}_{F(x_\pm)}\!\left[(\Lambda_{N_\wedge}-1)^2\right]\right)^{1/2}
    \left(\mathbb{E}_{F(x_\pm)}\!\left[R_\pm^2\,\mathbf{1}_{E_\pm}\right]\right)^{1/2}
= O\!\left(h^{3/2}\log^2(1/h)\right),
\end{equation*}
proving Eq.~\eqref{eq:girsanov_overshoot_shift}. Finally, writing $\mathbb{E}_F[R_\pm \mid \text{exit at } x_\pm] = \big(\mathbb{E}_F[R_\pm \mathbf{1}_{E_\pm}] + O(h^q)\big)/\mathbb{P}_F(\text{exit at } x_\pm)$ (with $q > 2$) and dividing by $\mathbb{P}_F(\text{exit at } x_\pm) = p_\pm(x_0) + O(\sqrt{h}) > 0$ yields Eq.~\eqref{eq:girsanov_cond_shift}.
\end{proof}

We now combine the Dynkin identity, the linear-drift overshoot expansion, and Lemma~\ref{lem:girsanov_bl} to obtain the two-term asymptotic expansion for nonlinear drifts. We state the result explicitly as a \emph{formal matched-asymptotic expansion}: the component lemmas are proved above, while the identification of the boundary-layer entrance law with the stationary renewal measure and the matching to the global discrete splitting probabilities at $O(h^{3/2})$ are carried out via matched boundary-layer asymptotics and validated numerically in Section~\ref{sec:numerics}.

\begin{asymptotic}[Two-term MFPT bias for nonlinear Langevin equations]\label{prop:nonlinear}
Let $F \in C^3([x_-, x_+])$ with constant diffusivity $D = \sigma^2/2 > 0$, and fix $x_0 \in (x_-, x_+)$ separated from $x_\pm$ by $O(1)$ distances. For a weak order-$p$ scheme driven by Gaussian increments ($p=1$ for Euler--Maruyama; $p=2$ for Heun) with step $h \to 0^+$ in the small boundary-drift regime $|\mu_\pm| = |F(x_\pm)|\sqrt{h}/\sigma \ll 1$, matched boundary-layer asymptotics yield
\begin{equation}\label{eq:master_drift}
\tau(h) - \tau(0) = c_{1/2}(x_0)\,\sqrt{h} + C_1(x_0)\,h + O(h^{3/2}),
\end{equation}
with explicit coefficients
\begin{align}
c_{1/2}(x_0) &= \beta\sigma\Big(p_+(x_0)\,|u'(x_+)| + p_-(x_0)\,|u'(x_-)|\Big), \label{eq:c_half} \\
C_1(x_0) &= \sum_{s \in \{+,-\}} p_s(x_0)\!\left[
\frac{s}{4} F(x_s)\,|u'(x_s)|
- \left(\beta^2 + \frac{1}{4}\right) D\,u''(x_s)
\right] \nonumber \\
&\quad + \beta^2 \sigma^2 \big(p_-(x_0)v'(x_-) - p_+(x_0)v'(x_+)\big)\big(|u'(x_+)| - |u'(x_-)|\big)
+ \delta_{p,1}\,C_{\mathrm{int}}(x_0), \label{eq:C1}
\end{align}
where $\beta = -\zeta(1/2)/\sqrt{2\pi}$, $u(x)$ solves $F u' + D u'' = -1$ with $u(x_\pm)=0$, $v(x) = p_+(x)$ solves $F v' + D v'' = 0$ with $v(x_-)=0$ and $v(x_+)=1$, $p_-(x) = 1 - v(x)$, $u''(x_\pm) = -\big(1 + F(x_\pm)u'(x_\pm)\big)/D$, and $C_{\mathrm{int}}(x_0)$ is the interior weak-error integral for $p=1$ (vanishing for $p \ge 2$).
\end{asymptotic}

\noindent\emph{Derivation.}
Starting from the Dynkin and backward-Kolmogorov identities in Eqs.~\eqref{eq:dynkin}--\eqref{eq:lhs_drift}, the derivation proceeds in three steps:

\noindent\textbf{Step 1 (Linear-drift expansion of $\E_\mu[R]$).}
For a unit-step Gaussian walk with small outward drift $\mu > 0$, Spitzer's identity gives the renewal overshoot expectation \cite{Siegmund1979,Siegmund1985}:
\begin{equation}\label{eq:siegmund_exact}
\E_\mu[R] = \frac{\mu}{2} + \frac{1}{2\mu} - \sum_{n=1}^\infty \frac{1}{\sqrt{n}}\Big[\varphi(\mu\sqrt{n}) - \mu\sqrt{n}\,\Phi(-\mu\sqrt{n})\Big],
\end{equation}
where $\Phi$ is the standard normal cumulative distribution function. Setting $g(u) := \varphi(u) - \varphi(0) - u\Phi(-u)$, splitting off $\varphi(0) = (2\pi)^{-1/2}$ via $\sum_{n=1}^N n^{-1/2} = 2\sqrt{N} + \zeta(1/2) + O(N^{-1/2})$, and evaluating the Mellin-transform poles of $\sum_{n=1}^\infty n^{-1/2} g(\mu\sqrt{n})$ at $s = 1$ ($\frac{2}{\mu}\int_0^\infty g(u)\dd u = \frac{1}{2\mu}$), $s = 0$ ($-\frac{1}{2}\zeta(0)\mu = +\frac{1}{4}\mu$ since $\zeta(0) = -\frac{1}{2}$), and $s = -1$ ($O(\mu^2)$) cancels $\frac{1}{2\mu}$ identically and recovers the classical linear drift expansion \cite{Siegmund1979,ChangPeres1997}
\begin{equation}\label{eq:overshoot_mu}
\E_\mu[R] = \beta + \frac{1}{4}\mu + O(\mu^2).
\end{equation}
For inward drift ($\mu < 0$), the ascending ladder height $H$ is defective with $\P_\mu(T_+ < \infty) = 1 - \sqrt{2}|\mu| + O(\mu^2) = e^{\sqrt{2}\mu + O(\mu^2)}$ (from Wald's identity and duality, $\P_\mu(T_+ = \infty) = |\mu|/\E_0[H] + O(\mu^2)$), and the conditional ladder-height expansion of Chang and Peres~\cite{ChangPeres1997} and Lotov~\cite{Lotov1996} shows that Eq.~\eqref{eq:overshoot_mu} extends analytically across $\mu = 0$ for both signs of $\mu$ as $\mu \to 0$.

\noindent\textbf{Step 2 (Girsanov invariance of the boundary layer to drift gradients $F'(x_\pm)$).}
Near $x_\pm$, rescale space via $x = x_\pm \mp \sigma\sqrt{h}\,y$ ($y > 0$). By Taylor's theorem, the step-normalized outward drift at distance $y$ from $x_\pm$ is $\mu(y) = \mu_\pm \mp F'(x_\pm)\,y\,h + O(h^{3/2}y^2)$ with $\mu_\pm = \pm F(x_\pm)\sqrt{h}/\sigma = O(\sqrt{h})$. By Lemma~\ref{lem:girsanov_bl}, replacing the variable-drift process $\P_F$ inside the boundary layer by the frozen-drift walk $\P_{\mu_\pm}$ changes the conditional overshoot moments only by $O(h^{3/2}\log^2(1/h))$---showing that spatial drift gradients $F'(x_\pm)$ do not generate any compensating $O(h)$ term in the local overshoot. Matching the boundary-layer entrance law to the stationary renewal distribution of the frozen-drift walk and applying Eq.~\eqref{eq:overshoot_mu} at $\mu = \mu_\pm$ therefore gives
\begin{equation}\label{eq:cond_moments}
\E[R_\pm \mid E_\pm] = \beta\,\sigma\sqrt{h} \pm \frac{1}{4}F(x_\pm)\,h + O(h^{3/2}),
\qquad
\E[R_\pm^2 \mid E_\pm] = \Big(\beta^2 + \frac{1}{4}\Big)\sigma^2 h + O(h^{3/2}).
\end{equation}
A priori, one might expect the local drift gradient $F'(x_\pm)$ to contribute at order $O(h)$, because the discrete boundary layer has spatial width $O(\sqrt{h})$ and a drift variation $\Delta F \sim F'(x_\pm)\sqrt{h}$ integrated over $O(1)$ boundary-layer steps produces an $O(h)$ perturbation in step increments. Lemma~\ref{lem:girsanov_bl} shows that this variation acts only as an $O(h)$ Radon--Nikodym change of measure on an overshoot $R_\pm$ whose amplitude is already $O(\sqrt{h})$, and therefore shifts the mean exit time only at order $O(h^{3/2})$.

\noindent\textbf{Step 3 (Splitting-probability shift via Dynkin's identity and assembly).}
Applying Dynkin's telescoping identity to the harmonic splitting function $v(x)$ ($Fv' + Dv'' = 0$, $v(x_-)=0$, $v(x_+)=1$) gives $\E_{x_0}[v(S_{N_h})] - v(x_0) = O(h)$. Expanding $v(x_\pm \pm R_\pm)$ at exit yields
\[
p_+(h)\Big(1 + v'(x_+)\,\E[R_+ \mid E_+] + O(h)\Big) - \big(1 - p_+(h)\big)\Big(v'(x_-)\,\E[R_- \mid E_-] + O(h)\Big) = p_+ + O(h).
\]
Substituting Eq.~\eqref{eq:cond_moments} and solving for $p_\pm(h)$ gives $p_\pm(h) = p_\pm \pm \beta\sigma\sqrt{h}\big(p_- v'(x_-) - p_+ v'(x_+)\big) + O(h)$. Combining this shift with Eqs.~\eqref{eq:dynkin}, \eqref{eq:lhs_drift}, and \eqref{eq:cond_moments} yields Eqs.~\eqref{eq:master_drift}--\eqref{eq:C1}. \hfill$\square$

\subsection{Shifted-Domain Equivalence and Universal \texorpdfstring{$h/4$}{h/4} Offset}\label{sec:shifted_domain}

Substituting $D u''(x_s) = -1 + s F(x_s)|u'(x_s)|$ into Eq.~\eqref{eq:C1}, the linear-drift term $\frac{s}{4}F(x_s)|u'(x_s)|$ cancels $-\frac{1}{4}D u''(x_s)$ identically, simplifying Eq.~\eqref{eq:C1} to
\begin{align}\label{eq:C1_simplified}
C_1(x_0) &= \frac{1}{4} - \beta^2 D\Big(p_+(x_0) u''(x_+) + p_-(x_0) u''(x_-)\Big) \nonumber \\
&\quad + \beta^2 \sigma^2 \big(p_-(x_0)v'(x_-) - p_+(x_0)v'(x_+)\big)\big(|u'(x_+)| - |u'(x_-)|\big)
+ \delta_{p,1}\,C_{\mathrm{int}}(x_0).
\end{align}
Remarkably, this three-term coefficient is algebraically equivalent to a second-order domain perturbation for \emph{arbitrary} nonlinear drift $F \in C^3([x_-, x_+])$.

\begin{proposition}[Shifted-domain equivalence and universal $h/4$ barrier-shift offset]\label{prop:shifted_domain}
Under the hypotheses of Asymptotic Expansion~\ref{prop:nonlinear} for a weak order-$p$ integrator driven by Gaussian increments ($p = 1$ for Euler--Maruyama; $p = 2$ for stochastic Heun), the two-term expansion~\eqref{eq:master_drift}--\eqref{eq:C1} is algebraically equivalent to the continuum MFPT on a domain widened by $\beta\sigma\sqrt{h}$ plus a universal $+\frac{h}{4}$ variance offset:
\begin{equation}\label{eq:shifted_domain_nonlinear}
\tau(h) = u\big(x_0;\, x_- - \beta\sigma\sqrt{h},\, x_+ + \beta\sigma\sqrt{h}\big) + \frac{h}{4} + \delta_{p,1}\,C_{\mathrm{int}}(x_0)\,h + O(h^{3/2}),
\end{equation}
where $u(x;\, a, b)$ solves $F(x)u'(x) + D u''(x) = -1$ with $u(a) = u(b) = 0$. Consequently, applying the inward barrier shift $\tilde{x}_\pm = x_\pm \mp \beta\sigma\sqrt{h}$ to a weak-order-2 scheme driven by Gaussian increments leaves a universal $O(h)$ residual independent of $F(x)$, $D$, and $x_0$:
\begin{equation}\label{eq:bs_universal}
\tau_{\mathrm{bs}}(h) = \tau(0) + \frac{h}{4} + O(h^{3/2}) \qquad (p \ge 2,\;\text{Gaussian increments}),
\end{equation}
so that the offset estimator $\tau_{\mathrm{bs}}(h) - h/4$ eliminates both the $O(\sqrt{h})$ and $O(h)$ bias simultaneously at zero extra computational cost.
\end{proposition}

\begin{proof}
Let $u_\delta(x) := u(x;\, x_- - \delta,\, x_+ + \delta)$ solve $F u_\delta' + D u_\delta'' = -1$ with $u_\delta(x_- - \delta) = u_\delta(x_+ + \delta) = 0$, and expand $u_\delta(x) = u(x) + \delta u_1(x) + \frac{1}{2}\delta^2 u_2(x) + O(\delta^3)$. Because the generator $L = F\partial_x + D\partial_{xx}$ is independent of $\delta$, both $u_1$ and $u_2$ are harmonic ($L u_1 = L u_2 = 0$), so $u_k(x_0) = p_+(x_0)u_k(x_+) + p_-(x_0)u_k(x_-)$ for $k = 1, 2$. Expanding the boundary conditions $u_\delta(x_s + s\delta) = 0$ ($s \in \{+,-\}$) gives $u_1(x_s) = -s u'(x_s) = |u'(x_s)|$, whence $u_1(x) = p_+(x)|u'(x_+)| + p_-(x)|u'(x_-)|$ and $u_1'(x) = (|u'(x_+)| - |u'(x_-)|)v'(x)$. At $O(\delta^2)$, $u_2(x_s) = -u''(x_s) - 2s u_1'(x_s)$. Setting $\delta = \beta\sigma\sqrt{h}$ (with $\frac{1}{2}\delta^2 = \beta^2 D h$) reproduces every $\beta$ and $\beta^2$ term of Eqs.~\eqref{eq:c_half} and~\eqref{eq:C1_simplified} identically, establishing Eq.~\eqref{eq:shifted_domain_nonlinear}. Replacing $(x_-, x_+)$ by the inward-shifted boundaries $\tilde{x}_\pm = x_\pm \mp \beta\sigma\sqrt{h}$ maps the effective widened domain $(\tilde{x}_- - \beta\sigma\sqrt{h},\, \tilde{x}_+ + \beta\sigma\sqrt{h})$ back onto $(x_-, x_+)$ identically, yielding Eq.~\eqref{eq:bs_universal} for $p \ge 2$.
\end{proof}

Two remarks clarify the physical origin and practical scope of Proposition~\ref{prop:shifted_domain}. First, the $+\frac{h}{4}$ offset originates directly from the Gaussian overshoot variance $\operatorname{Var}(R) = \frac{1}{4}\sigma^2 h$; for free Brownian motion on $[-L,L]$, Wald's identity~\eqref{eq:ratio} together with the asymptotic overshoot moments in Eq.~\eqref{eq:free_final} provides an immediate exact check: shifting the barrier inward to $L' = L - \beta\sigma\sqrt{h}$ gives $\tau_{\mathrm{bs}}(h) = [(L' + \beta\sigma\sqrt{h})^2 + \frac{1}{4}\sigma^2 h]/\sigma^2 = L^2/\sigma^2 + \frac{h}{4}$ up to exponentially small terms. Second, the constants $\beta$ and $\frac{1}{4}$ are specific to \emph{Gaussian} noise increments: weak-order-2 integrators that substitute discrete two- or three-point random variables for Gaussian increments match the low-order step moments in the interior but alter the discrete ladder-height distribution at the boundary, changing both the leading shift constant and the $O(h)$ overshoot-variance offset. Likewise, for Euler--Maruyama ($p = 1$), Eq.~\eqref{eq:shifted_domain_nonlinear} shows that the barrier-shifted residual is $\big(\frac{1}{4} + C_{\mathrm{int}}(x_0)\big)h$ rather than $h/4$.

\begin{remark}[Scope, boundary layers, and extensions]\label{rem:scope}
Asymptotic Expansion~\ref{prop:nonlinear} and Proposition~\ref{prop:shifted_domain} assume $F \in C^3([x_-, x_+])$ and an initial point $x_0$ separated from $x_\pm$ by $O(1)$ distances. Relaxing these hypotheses modifies the expansion in predictable ways:
\begin{enumerate}
  \item \emph{Piecewise-smooth drifts:} If $F$ or $F'$ has a jump discontinuity at an interior point $x_d \in (x_-, x_+)$ while remaining $C^3$ in neighborhoods of $x_\pm$, the local boundary-layer structure at $x_\pm$ is unaffected, so the leading $O(\sqrt{h})$ overshoot term in Eq.~\eqref{eq:master_drift} remains valid. However, repeated crossings of the kink at $x_d$ degrade the interior weak order of the Heun integrator from $O(h^2)$ to $O(h^{1/2})$ or $O(h)$, introducing an additional interior $O(h)$ contribution to $C_1(x_0)$.
  \item \emph{Initial condition inside the $O(\sqrt{h})$ boundary layer:} When $|x_0 - x_\pm| = a_0 \sigma\sqrt{h}$ with $a_0 = O(1)$, both $\tau(0) = O(\sqrt{h})$ and the exit-time discretization error $\tau(h) - \tau(0) = O(\sqrt{h})$ share the same order of magnitude, so the relative bias $\tau(h)/\tau(0) - 1$ becomes $O(1)$ rather than $O(\sqrt{h})$. In this regime, $p_\pm(x_0) = O(\sqrt{h})$ and one must replace the continuum splitting probabilities in Asymptotic Expansion~\ref{prop:nonlinear} with the exact discrete one-sided Wiener--Hopf exit probability and overshoot starting at scaled distance $a_0$.
  \item \emph{Multiplicative noise:} For a non-degenerate state-dependent diffusion $\mathrm{d}X_t = F(X_t)\,\mathrm{d}t + \sigma(X_t)\,\mathrm{d}W_t$ with $\sigma(x) > 0$ on $[x_-, x_+]$, the Lamperti transformation $Y_t = \int_{x_0}^{X_t} \mathrm{d}z/\sigma(z)$ maps the process onto an additive-noise Langevin equation with unit diffusivity and transformed drift $\tilde{F}(y) = F(x(y))/\sigma(x(y)) - \frac{1}{2}\sigma'(x(y))$. Applying Asymptotic Expansion~\ref{prop:nonlinear} in $Y$-coordinates yields the explicit $O(\sqrt{h})$ and $O(h)$ bias coefficients (with local boundary step scale $\sigma(x_\pm)\sqrt{h}$), provided the time-stepping scheme is applied in $Y$-coordinates or augmented with its standard It\^o--Stratonovich drift correction.
  \item \emph{Multidimensional domains:} In $\Omega \subset \mathbb{R}^d$ with smooth boundary $\partial\Omega$, projecting the local dynamics along the inward unit normal $\mathbf{n}(\mathbf{y})$ at $\mathbf{y} \in \partial\Omega$ recovers the leading $O(\sqrt{h})$ normal boundary-layer shift $\beta\sqrt{2 D_n(\mathbf{y}) h}$ with $D_n(\mathbf{y}) = \mathbf{n}(\mathbf{y})^\top \mathbf{D}(\mathbf{y})\mathbf{n}(\mathbf{y})$ (consistent with Gobet~\cite{Gobet2000}). Gobet and Menozzi~\cite{GobetMenozzi2010} showed that removing the $O(\sqrt{h})$ error via boundary shifting requires a normal shift incorporating the local mean curvature $\nabla\cdot\mathbf{n}(\mathbf{y})$ at order $O(h)$; we conjecture that in unshifted multidimensional grid thresholding, tangential diffusion along $\partial\Omega$ likewise introduces an explicit geometric contribution proportional to $\nabla\cdot\mathbf{n}(\mathbf{y})$ in $C_1(\mathbf{x}_0)$.
\end{enumerate}
\end{remark}

\subsection{Relation to Stopped-Diffusion Weak-Error Theory}\label{sec:weak_error_comp}

While the $O(\sqrt{h})$ weak convergence rate of killed time-stepping schemes is well established in numerical probability, existing analyses either treat constant-coefficient models, bound the $O(h)$ remainder implicitly, or modify the boundary rule to suppress the $O(\sqrt{h})$ error:
\begin{itemize}
  \item \emph{Constant-drift corrected diffusion approximations:} Siegmund~\cite{Siegmund1979,Siegmund1985} derived two-sided corrected diffusion approximations for $\E[N_h]$ in constant-drift Gaussian walks ($F(x) \equiv F_0$), and Broadie, Glasserman, and Kou~\cite{BroadieGlassermanKou1997} established the barrier shift $\beta\sigma\sqrt{h}$ for discretely monitored geometric Brownian motion. Proposition~\ref{prop:shifted_domain} extends this two-sided shifted-domain structure to arbitrary smooth nonlinear drifts $F(x)$ and identifies the universal $+\frac{h}{4}$ variance offset.
  \item \emph{Malliavin weak-error expansion of Gobet~\cite{Gobet2000}:} Extending weak-error expansions of unkilled diffusions~\cite{BallyTalay1996} via Malliavin integration by parts for killed Euler schemes on $\Omega \subset \mathbb{R}^d$, Gobet proved that the weak error expands as $c_{1/2}^{\mathrm{G}}\sqrt{h} + O(h)$, expressing $c_{1/2}^{\mathrm{G}}$ as a boundary integral weighted by the stationary overshoot $\beta\sigma$, while bounding the $O(h)$ remainder through Malliavin weight and occupation-time estimates rather than resolving its prefactor. In one dimension, $c_{1/2}(x_0)$ matches Gobet's leading term, whereas Eqs.~\eqref{eq:C1} and~\eqref{eq:C1_simplified} resolve the $O(h)$ coefficient of unshifted grid thresholding explicitly.
  \item \emph{Boundary-shifting, sub-step, and strong-convergence analyses:} Gobet and Menozzi~\cite{GobetMenozzi2004,GobetMenozzi2010} and Milstein and Tretyakov~\cite{MilsteinTretyakov2004} modified the discrete stopping boundary (via $\beta\sigma\sqrt{h}$ plus mean-curvature and drift shifts in $\mathbb{R}^d$) or the local transition law to eliminate the $O(\sqrt{h})$ bias and establish $O(h)$ weak-error bounds without identifying the explicit $O(h)$ prefactor of unshifted grid thresholding, while Buchmann~\cite{Buchmann2005} and Bernal and Acebr\'on~\cite{BernalAcebron2016} analyzed Brownian-bridge and local-linearization schemes. In a complementary strong-error direction, Bouchard, Geiss, and Gobet~\cite{BouchardGeissGobet2017} established general moment estimates and an $O(h^{1/2})$ $L_1$-convergence rate for the discrete Euler exit time.
\end{itemize}

\begin{table}[h]
\centering
\footnotesize
\setlength{\tabcolsep}{7pt}
\begin{tabular}{lccc}
\toprule
Property / Result & Gobet~\cite{Gobet2000} & Gobet--Menozzi~\cite{GobetMenozzi2004,GobetMenozzi2010} & Present work \\
\midrule
$O(\sqrt{h})$ weak error rate & $\checkmark$ & $\checkmark$ & $\checkmark$ \\
Explicit $c_{1/2}$ coefficient & $\checkmark$ & $\checkmark$ & $\checkmark$ \\
Explicit $O(h)$ coefficient $C_1$ & --- & --- & $\checkmark$ (1D) \\
Nonlinear drift $F(x)$ & $\checkmark$ & $\checkmark$ & $\checkmark$ \\
\bottomrule
\end{tabular}
\caption{Comparison of stopped-diffusion weak-error results for discretely monitored SDEs. Previous analyses prove the $O(\sqrt{h})$ expansion structure, evaluate the leading prefactor, or construct boundary-shifted schemes with proven $O(h)$ weak-error bounds, whereas Eqs.~\eqref{eq:C1}--\eqref{eq:bs_universal} give the explicit one-dimensional $O(h)$ bias coefficient of unshifted grid thresholding and show that the remaining $O(h)$ residual after domain shifting is universally $+h/4$ for weak-order-2 schemes driven by Gaussian increments.}
\label{tab:lit_comparison}
\end{table}

While the effective $\beta\sigma\sqrt{h}$ domain shift is classical~\cite{Siegmund1979,BroadieGlassermanKou1997,GobetMenozzi2010}, Eqs.~\eqref{eq:C1}--\eqref{eq:bs_universal} provide, to the best of our knowledge, the first explicit one-dimensional representation of the $O(h)$ bias coefficient for unshifted schemes under general nonlinear drift $F(x)$ and establish that the remaining $O(h)$ error after domain shifting is universally $+\frac{h}{4}$ for weak-order-2 integrators driven by Gaussian increments.

%======================================================================
\section{Numerical verification and boundary-correction performance}\label{sec:numerics}
%======================================================================

\subsection{Free Brownian motion and systematic quadrature convergence studies}

Table~\ref{tab:free} compares exact Nystr\"om Fredholm quadrature ($600$ Gauss--Legendre nodes) and OpenMP-parallelized Fortran~90 Monte Carlo simulations ($N = 5 \times 10^8$ trajectories per $a$) of the midpoint-started unit-step Gaussian walk on $[-a,a]$ ($h/\tau(0) = 1/a^2$) against Eq.~\eqref{eq:free_final} and the linear-density heuristic. At $a = 1$ ($h = \tau(0)$), the two-term asymptotic formula ($2.75461$) differs from the exact value ($2.78293$) by $\approx 1\%$ due to the exponential remainder $\mathcal{R}(h)$; already by $a = 3$ ($1.45389$ vs $1.45396$) the discrepancy is below $10^{-4}$, and for $a \ge 5$ all simulated moments agree with $\beta = 0.58260$ and $\beta^2 + \frac{1}{4} = 0.58942$ within statistical error.

\begin{table}[h]
\centering
\setlength{\tabcolsep}{4.5pt}
\footnotesize
\begin{tabular}{rcccccc}
\toprule
$a$ & $\E[R]$ & $\E[R^2]$ & $\tau(h)/\tau(0)$ & 1st order & 2nd order & Heuristic \\
    & (exact / MC) & (exact / MC) & (exact / MC) & $1 + 2\beta/a$ & Eq.~\eqref{eq:free_final} & $1 + \frac{4}{3}\sqrt{2/\pi}/a$ \\
\midrule
1  & $0.59181$ & $0.59930$ & $2.78293$ & $2.16519$ & $2.75461$ & $2.06385$ \\
2  & $0.58145$ & $0.58769$ & $1.72838$ & $1.58260$ & $1.72995$ & $1.53192$ \\
3  & $0.58268$ & $0.58957$ & $1.45396$ & $1.38840$ & $1.45389$ & $1.35462$ \\
5  & $0.58259(2)$ & $0.58940(4)$ & $1.25662(4)$ & $1.23304$ & $1.25662$ & $1.21277$ \\
10 & $0.58264(2)$ & $0.58950(4)$ & $1.12241(4)$ & $1.11652$ & $1.12241$ & $1.10638$ \\
20 & $0.58258(2)$ & $0.58940(4)$ & $1.05973(4)$ & $1.05826$ & $1.05973$ & $1.05319$ \\
\midrule
Theory & $\beta = 0.58260$ & $\beta^2 + \frac{1}{4} = 0.58942$ & & & & \\
\bottomrule
\end{tabular}
\caption{Exact Fredholm quadrature ($a = 1, 2, 3$) and OpenMP Fortran~90 Monte Carlo benchmarks ($a = 5, 10, 20$; $N = 5 \times 10^8$ trajectories per $a$) for midpoint-started free Brownian motion on $[-a,a]$ compared with Eq.~\eqref{eq:free_final} and the linear-density heuristic. Numbers in parentheses denote one standard error in the last digit.}
\label{tab:free}
\end{table}

To validate both the $O(\sqrt{h})$ and $O(h)$ terms of Asymptotic Expansion~\ref{prop:nonlinear} with quantified numerical uncertainties---including in independent regimes not used to motivate the asymptotics---we solved the discrete-time Fredholm integral equation $\E_x[N_h] = 1 + \int_{x_-}^{x_+} p_h(x,y)\,\E_y[N_h]\dd y$ using Gauss--Legendre quadrature ($N_{\mathrm{GL}} \in \{800, 1200, 1600\}$). For Heun's integrator ($p = 2$, $C_{\mathrm{int}} = 0$), the single-step map $y = \Psi_h(x,\xi) = x + \frac{1}{2}[F(x) + F(x + F(x)h + \sigma\sqrt{h}\xi)]h + \sigma\sqrt{h}\xi$ has a non-Gaussian transition density $p_h(x,y) = \varphi(\xi(x,y))/|\partial_\xi \Psi_h(x,\xi(x,y))|$, which we evaluate at each quadrature pair $(x_i, y_j)$ by inverting $y_j = \Psi_h(x_i,\xi)$ via four damped Newton iterations (for Euler--Maruyama, $p_h(x,y)$ is directly Gaussian). Richardson extrapolation is performed by fitting $[\tau(h) - \tau(0) - c_{1/2}\sqrt{h}]/h = C_1 + c_{3/2}\sqrt{h} + c_2 h$ in powers of $\sqrt{h}$:
\begin{enumerate}
\item \emph{Constant drift} ($F(x) = 1.0$, $D = 0.5$, $(x_-,x_+) = (0,2)$, $x_0 = 1.0$): $\tau(0) = \tanh(1) = 0.761594$, leading $\sqrt{h}$ slope $c_{1/2}(x_0) = 0.688378$, and Eq.~\eqref{eq:C1} predicts $C_1(x_0) = 0.283984$. Richardson extrapolation across $h \in \{0.004, 0.002, 0.001\}$ yields $0.28398(2)$ (quadrature grid uncertainty $< 10^{-6}$), with log-log remainder exponent $1.45$ across $h \in [0.001, 0.008]$.
\item \emph{Ornstein--Uhlenbeck drift} ($F(x) = -\gamma x$, $F'(x) = -\gamma \neq 0$, $(x_-,x_+) = (0,2)$, $x_0 = 0.8$): for $\gamma = 1.5$, Eq.~\eqref{eq:C1} predicts $C_1(x_0) = 0.572415$ and Richardson extrapolation yields $0.57243(2)$ (remainder exponent $1.59$); for $\gamma = 3.0$, Eq.~\eqref{eq:C1} predicts $C_1(x_0) = 0.588123$ and Richardson extrapolation across $h \in \{0.002, 0.001, 0.0005\}$ yields $0.58815(5)$ (with the uncertainty estimated from the shift relative to $h \in \{0.004, 0.002, 0.001\}$; remainder exponent $1.60$), confirming the absence of local $F'(x_\pm)$ contributions to $C_1(x_0)$.
\item \emph{Independent periodic washboard potential} ($F(x) = 1.5 - 2.5\sin(2\pi x)$, $D = 0.3$, $(x_-,x_+) = (0,1.5)$, $x_0 = 0.5$, where $F'(x_-) = -5\pi$ and $F'(x_+) = +5\pi \neq 0$): $\tau(0) = 0.771124$, leading $\sqrt{h}$ slope $c_{1/2}(x_0) = 1.037980$, and Eq.~\eqref{eq:C1} predicts $C_1(x_0) = -0.113057$. Here the $O(h^{3/2})$ and $O(h^2)$ remainder terms cross zero near $h \approx 10^{-3}$ (with $|\tau(h) - \tau_{\mathrm{2nd}}(h)|$ decreasing from $1.04 \times 10^{-3}$ at $h = 0.008$ to $1.15 \times 10^{-6}$ at $h = 0.001$); refining the step sequence to $h \in \{0.001, 0.0005, 0.00025\}$ ($N_{\mathrm{GL}} = 1600$) yields Richardson-extrapolated $C_1 = -0.1134(3)$, confirming Eq.~\eqref{eq:C1} under strong periodic forcing.
\item \emph{Euler--Maruyama ($p = 1$) interior weak-error term:} Repeating the OU benchmark ($\gamma = 1.5$) with the Euler--Maruyama scheme ($p = 1$), solving $F C_{\mathrm{int}}' + D C_{\mathrm{int}}'' = -\mathcal{A}_1 u$ gives $C_{\mathrm{int}}(x_0) = -0.800611$ and $C_1^{\mathrm{EM}}(x_0) = 0.572415 - 0.800611 = -0.228196$, in agreement with Richardson extrapolation of the Gaussian-kernel Fredholm solver ($-0.22821(2)$).
\end{enumerate}

\subsection{Nonlinear Langevin Benchmarks and Remainder Scaling}

\begin{figure}[t]
  \centering
  \includegraphics[width=0.72\textwidth]{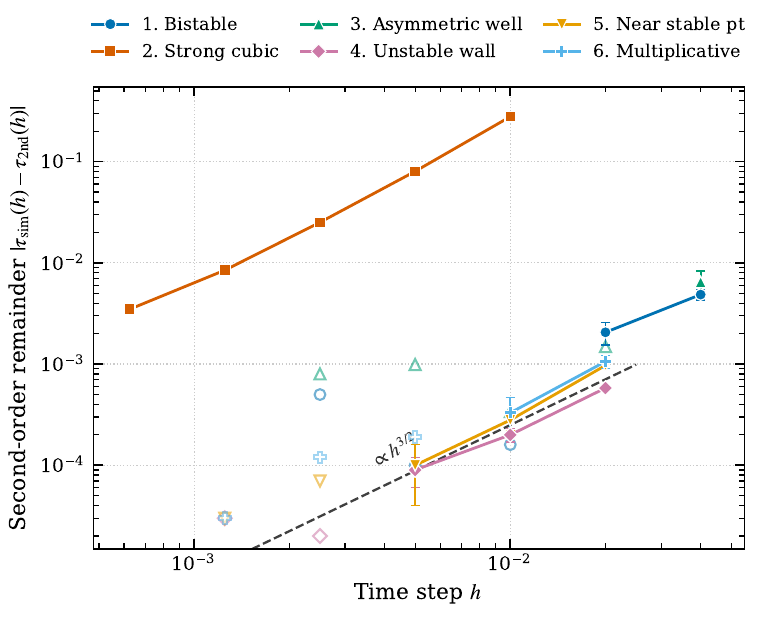}
\caption{Absolute remainder $|\tau_{\mathrm{sim}}(h) - \tau_{\mathrm{2nd}}(h)|$ between raw Heun simulations ($N = 10^8$) and the second-order asymptotic prediction $\tau_{\mathrm{2nd}}(h) = \tau(0) + c_{1/2}\sqrt{h} + C_1 h$ across all six nonlinear models. Filled symbols with $1\sigma$ error bars denote points above the Monte Carlo resolution threshold ($|\text{error}| > 1.5\,\mathrm{SE}$); open symbols indicate points limited by statistical sampling noise. Weighted least-squares log-log fits over the resolved regimes ($|\text{error}| > 2\,\mathrm{SE}$, $N_{\mathrm{pts}} \ge 3$) yield exponents $1.44 \pm 0.09$ for Model~4 (consistent with the dashed $O(h^{3/2})$ reference line) and $1.76 \pm 0.03$ for Model~2 as pre-asymptotic $O(h^2)$ terms subside, while Models~1, 3, 5, and~6 reach the Monte Carlo statistical noise floor ($\sim 10^{-4}$) for $h \le 0.01$.}
  \label{fig:remainder}
\end{figure}

\begin{figure}[t]
  \centering
  \includegraphics[width=\textwidth]{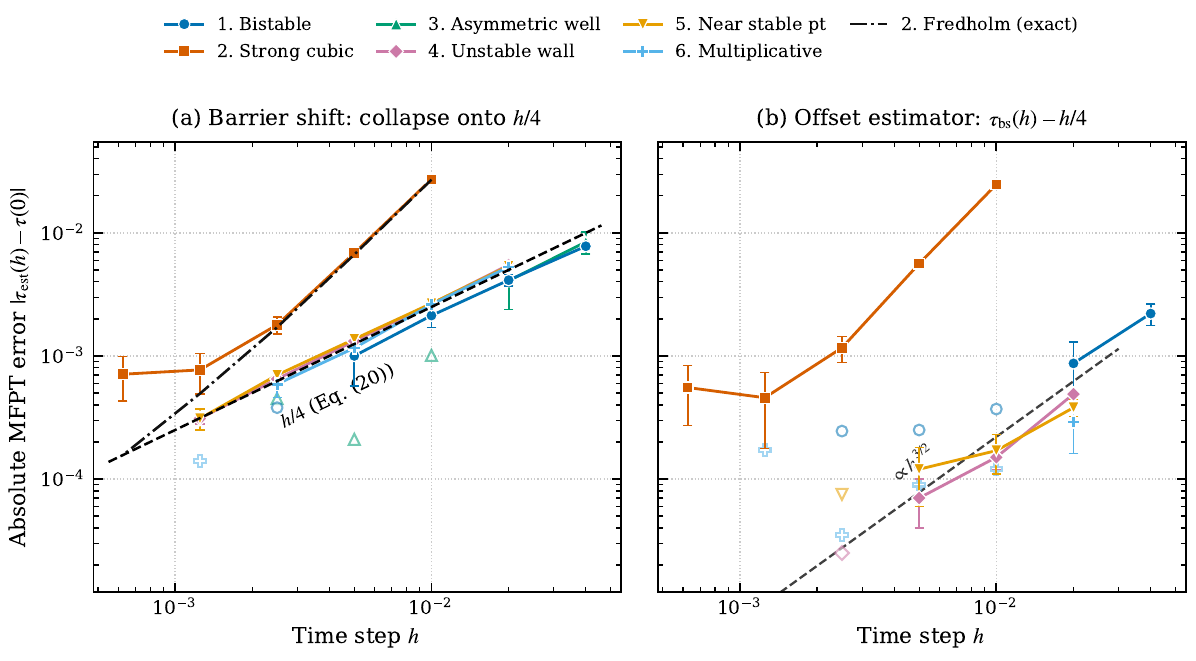}
  \caption{Universal $h/4$ barrier-shift offset and zero-cost $O(h^{3/2})$ correction across all six nonlinear models ($N = 10^8$). Left: residual MFPT error $|\tau_{\mathrm{bs}}(h) - \tau(0)|$ of the Broadie--Glasserman--Kou barrier shift $x_\pm \mp \beta\sigma\sqrt{h}$ compared against the exact parameter-free $h/4$ prediction of Eq.~\eqref{eq:bs_universal} (dashed line); the dash-dotted curve shows the noise-free Fredholm quadrature for Model~2, confirming the crossover from pre-asymptotic $O(h^2)$ decay toward $h/4$. Right: residual error $|\tau_{\mathrm{bs}}(h) - h/4 - \tau(0)|$ after subtracting $h/4$, plotted on the same vertical scale with a dashed $O(h^{3/2})$ reference line. Filled markers denote points above the Monte Carlo resolution threshold ($|\text{error}| > 1.5\,\mathrm{SE}$); open markers denote noise-floor-limited points. (Model~3 omitted from panel~(b) due to its large standard error)}
  \label{fig:corrections}
\end{figure}

Because the $O(h^{3/2})$ remainder in Asymptotic Expansion~\ref{prop:nonlinear} is derived via formal matched asymptotics, we benchmark Eqs.~\eqref{eq:c_half}--\eqref{eq:C1} against high-statistics Monte Carlo simulations ($N = 10^8\text{--}5\times 10^8$ trajectories) across six one-dimensional Langevin systems to validate the predicted coefficients and examine the residual scaling empirically. Each model is integrated using an OpenMP-parallelized Fortran~90 implementation of the stochastic Heun scheme coupled to independent per-thread \textsc{pcg64-dxsm} streams:
\begin{enumerate}
  \item \textbf{Bistable uphill escape} ($F(x) = x - x^3$, $D = 0.10$, $(x_-, x_+) = (-5.0, -0.50)$, $x_0 = -1.0$), representing thermally activated climb toward an absorbing threshold inside a quartic double well (note that although $F(x_-) = 120$ is large, $p_-(x_0) < 10^{-15}$ renders exits at $x_-$ negligible, while $|\mu_+| \approx 0.17$ at $h = 0.04$);
  \item \textbf{Strong cubic restoring drift} ($F(x) = 4x - 2x^3$, $D = 0.25$, $(x_-, x_+) = (-0.5, 1.80)$, $x_0 = 0.5$), featuring a steep inward drift $F(x_+) = -4.46$ at the upper exit boundary ($|\mu_+| \approx 0.63$ at $h = 0.01$, shrinking to $0.16$ at $h = 6.25 \times 10^{-4}$);
  \item \textbf{Asymmetric double well} ($F(x) = x - x^3 + 0.15$, $D = 0.08$, $(x_-, x_+) = (-2.5, 0.00)$, $x_0 = -0.9$), where tilted basin depths break exit symmetry;
  \item \textbf{Exit through an unstable fixed point} ($F(x) = 2x$, $D = 0.20$, $(x_-, x_+) = (-1.2, 0.00)$, $x_0 = -0.5$), where $F(x_+) = 0$ and the drift accelerates trajectories toward the boundary;
  \item \textbf{Exit near a stable fixed point} ($F(x) = -2x$, $D = 0.20$, $(x_-, x_+) = (-1.2, 0.05)$, $x_0 = -0.4$), where a linear restoring force opposes escape near $x_+ = 0.05$.
\end{enumerate}
Finally, to test the extension beyond additive-noise Langevin equations (Remark~\ref{rem:scope}, item~3), Model~6 considers a multiplicative-noise logistic SDE, $\mathrm{d}X_t = r X_t(1 - \kappa X_t)\,\mathrm{d}t + \sigma_0 X_t\,\mathrm{d}W_t$ with $r = 1.2$, $\kappa = 0.5$, and $\sigma_0 = 0.6$ ($D = \sigma_0^2/2 = 0.18$) on $(x_-, x_+) = (0.3, 1.8)$ starting at $X_0 = 1.0$. Because the stochastic Heun scheme achieves weak order $2$ only for additive noise (for multiplicative noise it converges to the Stratonovich SDE with weak order $1$ unless augmented), we integrate Model~6 in the Lamperti coordinate $Y_t = \ln X_t$, which satisfies the additive-noise Langevin equation $\mathrm{d}Y_t = \tilde{F}(Y_t)\,\mathrm{d}t + \sigma_0\,\mathrm{d}W_t$ with $\tilde{F}(y) = 1.02 - 0.6\,e^y$ on $(\ln 0.3, \ln 1.8)$ starting at $y_0 = 0$. As shown in Table~\ref{tab:nonlinear} and Figs.~\ref{fig:remainder}--\ref{fig:corrections}, the transformed coefficients $(c_{1/2}, C_1) = (1.30547, 0.66647)$ capture the discrete-time Heun bias to within $O(h^{3/2})$, while both the barrier shift and the Brownian-bridge exit test reduce the error to $O(h)$.

For each model, we sweep a geometric sequence of five time steps $h_k = h_0 2^{-k}$ ($k = 0,\dots,4$) and compare four exit estimators: (i)~uncorrected (\emph{raw}) Heun stopping on $X_{n+1} \notin (x_-, x_+)$, (ii)~the Broadie--Glasserman--Kou \emph{barrier shift} $x_\pm \mp \beta\sigma\sqrt{h}$ with $\beta = -\zeta(1/2)/\sqrt{2\pi} \approx 0.58260$, (iii)~the zero-cost $h/4$-offset barrier-shift estimator $\tau_{\mathrm{bs}}(h) - h/4$ of Eq.~\eqref{eq:bs_universal}, and (iv)~the standard drift-free \emph{Brownian-bridge} sub-step exit test of Ref.~\cite{Mannella1999} (which rejects interior steps $X_{n+1} \in (x_-, x_+)$ with probability $\exp[-2(x_+ - X_n)(x_+ - X_{n+1})/(\sigma^2 h)]$, and symmetrically at $x_-$). Continuum values $\tau(0)$ and derivatives $u'(x_\pm), v'(x_\pm)$ are obtained by solving the stationary boundary-value problems on a fourth-order finite-difference grid ($N_{\mathrm{grid}} = 4 \times 10^4$, verified against $N_{\mathrm{grid}} = 10^5$ to better than $3 \times 10^{-8}$; for Model~2, $\tau(0) = 3.1982618$).

\begin{table}[t]
  \centering
  \caption{Mean first-passage times for six nonlinear Langevin models integrated with the stochastic Heun scheme over five geometric step sizes ($N = 10^8$ trajectories per entry; parenthetical digits denote $1\sigma$ Monte Carlo standard errors in the last printed decimal places). Continuum $\tau(0)$, first-order $\tau_{\mathrm{1st}}(h) = \tau(0) + c_{1/2}\sqrt{h}$, and second-order $\tau_{\mathrm{2nd}}(h) = \tau_{\mathrm{1st}}(h) + C_1 h$ predictions are compared against raw Heun, barrier-shifted Heun ($x_\pm \mp \beta\sigma\sqrt{h}$), the $h/4$-offset barrier-shift estimator ($\mathrm{BS} - h/4$), and drift-free Brownian-bridge-corrected Heun.}
  \label{tab:nonlinear}
  \scriptsize
  \setlength{\tabcolsep}{2.8pt}
  \begin{tabular}{lcccccccc}
    \toprule
    Model & $h$ & $\tau(0)$ & $\tau_{\mathrm{1st}}(h)$ & $\tau_{\mathrm{2nd}}(h)$ & Raw Heun & Barrier Shift & BS $- h/4$ & Br.\ Bridge \\
    \midrule
    1.\ Bistable      & $0.04000$ & $4.55066$ & $5.73860$ & $5.87824$ & $5.88310(56)$ & $4.55846(43)$ & $4.54846(43)$ & $4.55897(43)$ \\
      ($x_+=-0.5$)                         & $0.02000$ & $4.55066$ & $5.39066$ & $5.46048$ & $5.46254(52)$ & $4.55480(43)$ & $4.54980(43)$ & $4.55417(43)$ \\
                               & $0.01000$ & $4.55066$ & $5.14463$ & $5.17954$ & $5.17970(49)$ & $4.55280(43)$ & $4.55030(43)$ & $4.55243(43)$ \\
                               & $0.00500$ & $4.55066$ & $4.97066$ & $4.98812$ & $4.98822(48)$ & $4.55167(43)$ & $4.55042(43)$ & $4.55160(43)$ \\
                               & $0.00250$ & $4.55066$ & $4.84765$ & $4.85637$ & $4.85588(46)$ & $4.55105(43)$ & $4.55043(43)$ & $4.55036(43)$ \\
\midrule
2.\ Strong cubic & $0.01000$ & $3.19826$ & $4.76413$ & $5.34471$ & $5.62297(52)$ & $3.22553(28)$ & $3.22303(28)$ & $3.21580(28)$ \\
   ($F_+=-4.46$)                             & $0.00500$ & $3.19826$ & $4.30550$ & $4.59579$ & $4.67596(43)$ & $3.20513(28)$ & $3.20388(28)$ & $3.20335(28)$ \\
                               & $0.00250$ & $3.19826$ & $3.98120$ & $4.12634$ & $4.15148(37)$ & $3.20005(28)$ & $3.19943(28)$ & $3.19994(28)$ \\
                               & $0.00125$ & $3.19826$ & $3.75188$ & $3.82445$ & $3.83291(34)$ & $3.19903(28)$ & $3.19872(28)$ & $3.19837(28)$ \\
                               & $0.00063$ & $3.19826$ & $3.58973$ & $3.62602$ & $3.62949(32)$ & $3.19897(28)$ & $3.19881(28)$ & $3.19863(28)$ \\
\midrule
3.\ Asymmetric well            & $0.04000$ & $18.07921$ & $19.57858$ & $19.53664$ & $19.53017(184)$ & $18.08766(171)$ & $18.07766(171)$ & $18.08951(171)$ \\
                               & $0.02000$ & $18.07921$ & $19.13942$ & $19.11846$ & $19.11696(180)$ & $18.08329(171)$ & $18.07829(171)$ & $18.08233(171)$ \\
                               & $0.01000$ & $18.07921$ & $18.82889$ & $18.81841$ & $18.81807(178)$ & $18.08022(171)$ & $18.07772(171)$ & $18.08120(171)$ \\
                               & $0.00500$ & $18.07921$ & $18.60932$ & $18.60407$ & $18.60506(176)$ & $18.07942(171)$ & $18.07817(171)$ & $18.08110(171)$ \\
                               & $0.00250$ & $18.07921$ & $18.45405$ & $18.45143$ & $18.45063(174)$ & $18.07876(171)$ & $18.07814(171)$ & $18.07922(171)$ \\
\midrule
4.\ Unstable wall    & $0.02000$ & $0.42234$ & $0.45774$ & $0.46144$ & $0.46202(3)$ & $0.42783(3)$ & $0.42283(3)$ & $0.42866(3)$ \\
  ($x_b=0$)                              & $0.01000$ & $0.42234$ & $0.44737$ & $0.44922$ & $0.44942(3)$ & $0.42499(3)$ & $0.42249(3)$ & $0.42529(3)$ \\
                               & $0.00500$ & $0.42234$ & $0.44004$ & $0.44096$ & $0.44105(3)$ & $0.42366(3)$ & $0.42241(3)$ & $0.42377(3)$ \\
                               & $0.00250$ & $0.42234$ & $0.43486$ & $0.43532$ & $0.43534(3)$ & $0.42299(3)$ & $0.42237(3)$ & $0.42301(3)$ \\
                               & $0.00125$ & $0.42234$ & $0.43119$ & $0.43142$ & $0.43145(3)$ & $0.42265(3)$ & $0.42234(3)$ & $0.42267(3)$ \\
\midrule
5.\ Near stable pt    & $0.02000$ & $0.63266$ & $0.75226$ & $0.76534$ & $0.76631(6)$ & $0.63804(6)$ & $0.63304(6)$ & $0.63792(6)$ \\
    ($x^*=0$)                           & $0.01000$ & $0.63266$ & $0.71723$ & $0.72377$ & $0.72405(6)$ & $0.63533(6)$ & $0.63283(6)$ & $0.63530(6)$ \\
                               & $0.00500$ & $0.63266$ & $0.69246$ & $0.69573$ & $0.69583(6)$ & $0.63403(6)$ & $0.63278(6)$ & $0.63387(6)$ \\
                               & $0.00250$ & $0.63266$ & $0.67494$ & $0.67658$ & $0.67665(6)$ & $0.63336(6)$ & $0.63274(6)$ & $0.63330(6)$ \\
                               & $0.00125$ & $0.63266$ & $0.66256$ & $0.66338$ & $0.66335(6)$ & $0.63297(6)$ & $0.63266(6)$ & $0.63285(6)$ \\
\midrule
6.\ Multiplicative  & $0.02000$ & $1.40351$ & $1.58813$ & $1.60146$ & $1.60252(15)$ & $1.40880(13)$ & $1.40380(13)$ & $1.40850(13)$ \\
  (Lamperti)                              & $0.01000$ & $1.40351$ & $1.53405$ & $1.54072$ & $1.54110(14)$ & $1.40613(13)$ & $1.40363(13)$ & $1.40611(13)$ \\
                               & $0.00500$ & $1.40351$ & $1.49582$ & $1.49915$ & $1.49946(14)$ & $1.40467(13)$ & $1.40342(13)$ & $1.40459(13)$ \\
                               & $0.00250$ & $1.40351$ & $1.46878$ & $1.47045$ & $1.47064(14)$ & $1.40410(13)$ & $1.40348(13)$ & $1.40409(13)$ \\
                               & $0.00125$ & $1.40351$ & $1.44966$ & $1.45050$ & $1.45052(14)$ & $1.40365(13)$ & $1.40334(13)$ & $1.40383(13)$ \\
    \bottomrule
  \end{tabular}
\end{table}

Table~\ref{tab:nonlinear} and Figs.~\ref{fig:remainder}--\ref{fig:corrections} demonstrate three quantitative features of the discretization error:
\begin{itemize}
  \item \textbf{Dominance of the $O(\sqrt{h})$ overshoot and accuracy of $\tau_{\mathrm{2nd}}(h)$:} Even at fine time steps ($h \sim 10^{-3}$), raw Heun integration overestimates $\tau(0)$ by $2\%$--$15\%$ due to the $c_{1/2}\sqrt{h}$ boundary-layer term. Incorporating the $O(h)$ drift-curvature and splitting-probability shift $C_1 h$ reduces the discrepancy between theory and simulation by one to two orders of magnitude across all six models. In Models~1, 3, 4, 5, and~6, $\tau_{\mathrm{2nd}}(h)$ agrees with the raw Heun estimates to within $\sim 10^{-4}\text{--}10^{-3}$ across the entire step-size sweep, reaching the Monte Carlo noise floor ($\sim 3\text{--}18 \times 10^{-4}$) for $h \le 0.01$.
  \item \textbf{Remainder scaling across models:} Figure~\ref{fig:remainder} displays the second-order remainder $|\tau_{\mathrm{sim}}(h) - \tau_{\mathrm{2nd}}(h)|$ on logarithmic axes. In Model~4 (where the small standard error $\approx 2.7 \times 10^{-5}$ resolves the remainder over multiple step sizes), a weighted log-log fit yields a scaling exponent of $1.44 \pm 0.09$, consistent with the predicted $O(h^{3/2})$ remainder scaling. In Model~2, the large boundary drift ($|\mu_+| \approx 0.63$ at $h = 0.01$) places the coarser steps outside the asymptotic $|\mu_+| \ll 1$ regime, producing a larger pre-asymptotic remainder ($0.27827$ at $h = 0.01$ decreasing to $0.00347$ at $h = 6.25 \times 10^{-4}$) with an effective slope of $1.76 \pm 0.03$ as $O(h^2)$ corrections subside, while Models~1, 3, 5, and~6 reach the statistical noise floor for $h \le 0.01$.
  \item \textbf{Reduction to $O(h^1)$ weak convergence and verification of the universal $h/4$ residual:} As shown in Table~\ref{tab:nonlinear}, both the static barrier shift ($x_\pm \mp \beta\sigma\sqrt{h}$) and the drift-free Brownian-bridge sub-step exit test remove the $O(\sqrt{h})$ bias and achieve comparable accuracy across all six test cases. Weighted log-log regressions of $|\tau_{\mathrm{corr}}(h) - \tau(0)|$ over the statistically resolved regimes ($|\text{error}| > 2\,\mathrm{SE}$) confirm linear $O(h^1)$ convergence: we obtain slopes of $0.94 \pm 0.02$ (barrier shift) and $1.14 \pm 0.06$ (Brownian bridge) for Model~1, $1.03 \pm 0.01$ and $1.08 \pm 0.01$ for Model~4, and $1.00 \pm 0.01$ and $1.03 \pm 0.03$ for Model~5 (the slight difference in pre-asymptotic slope between the two methods reflects the different $O(h)$ drift residuals left by a static shift versus a drift-free bridge). Figure~\ref{fig:corrections} illustrates both the universal collapse of the barrier-shifted residuals onto the exact $h/4$ line of Eq.~\eqref{eq:bs_universal} (left) and the order-of-magnitude error reduction achieved by subtracting $h/4$ (right). Across all statistically resolved models (Models~1, 2, 4, 5, and~6; in Model~3 the $1.7\times 10^{-3}$ standard error exceeds $h/4$ at the finest step and is therefore uninformative), the finest-step barrier-shift residuals match $h/4$: at $h = 1.25\times 10^{-3}$, Models~4 and~5 both give $\tau_{\mathrm{bs}}(h) - \tau(0) = 0.248\,h$ (within $0.1\sigma$ of $0.25\,h$), while Models~1 and~6 agree with $h/4$ within $0.5\sigma$ and $1.3\sigma$, respectively. As shown in the $\mathrm{BS} - h/4$ column of Table~\ref{tab:nonlinear} and Fig.~\ref{fig:corrections}(b), subtracting $h/4$ at $h = 0.02$ drops the residual error by more than an order of magnitude---from $+5.49\times 10^{-3}$ to $+4.9\times 10^{-4}$ in Model~4 ($11\times$), from $+5.38\times 10^{-3}$ to $+3.8\times 10^{-4}$ in Model~5 ($14\times$), and from $+5.29\times 10^{-3}$ to $+2.9\times 10^{-4}$ in Model~6 ($18\times$)---and reaches the continuum value $\tau(0)$ within $0.1\sigma$ by $h = 1.25\times 10^{-3}$ in Models~4 and~5. In Model~2, where $|\mu_+| \approx 0.63$ at $h = 0.01$, exact noise-free Fredholm quadrature of the barrier-shifted Heun scheme yields $\tau_{\mathrm{bs}}(h) - \tau(0) = 0.02707,\; 0.00667,\; 0.00170,\; 4.74\times 10^{-4},\; 1.56\times 10^{-4}$ across $h \in [0.01, 6.25\times 10^{-4}]$ (dotted curve in Fig.~\ref{fig:corrections}(a)), showing pre-asymptotic $O(h^2)$ decay at coarse steps approaching the universal $h/4$ asymptote at fine steps ($1.56\times 10^{-4} = 0.250\,h$ at $h = 6.25\times 10^{-4}$); the apparent plateau in the last two Monte Carlo entries ($3.19903(28)$ and $3.19897(28)$ versus $\tau(0) = 3.1982618$) is a $+2.0\sigma$ statistical fluctuation inside the $2.8\times 10^{-4}$ noise floor.
\end{itemize}

\paragraph{Limitations.}
The derivation of Asymptotic Expansion~\ref{prop:nonlinear} combines rigorous renewal-theoretic and localized Girsanov estimates with a formal matched-asymptotic identification of the boundary-layer entrance law and the $O(\sqrt{h})$ splitting-probability correction. While the algebraic equivalence between Eq.~\eqref{eq:master_drift} and the shifted-domain representation~\eqref{eq:shifted_domain_nonlinear} is exact, the resulting universal $h/4$ barrier-shift corollary~\eqref{eq:bs_universal} inherits the formal matched-asymptotic status and hypotheses ($|\mu_\pm| \ll 1$, Gaussian increments, and $x_0$ separated from $x_\pm$ by $O(1)$ distances) of Asymptotic Expansion~\ref{prop:nonlinear}. Although both high-precision Fredholm quadrature and $N = 10^8\text{--}5\times 10^8$ Monte Carlo simulations across six nonlinear models support the resulting two-term expansion and $O(h^{3/2})$ remainder scaling, a fully rigorous probabilistic proof of the $O(h^{3/2})$ global remainder---as well as an explicit evaluation of the $O(h)$ bias coefficient in multidimensional curved domains---remains open.

%======================================================================
\section{Conclusions}\label{sec:conclusions}
%======================================================================

Discrete-time monitoring of a continuous diffusion produces a systematic $O(\sqrt{h})$ positive bias in the mean exit time whenever naive grid-point thresholding is applied. The contributions of this work are threefold: (i)~a self-contained Mellin-transform derivation of the first two overshoot moments for midpoint-started free Brownian motion and a discrete Milne boundary-layer calculation explaining why continuum linear-density heuristics underestimate the leading prefactor by $\approx 9\%$; (ii)~the matched-asymptotic expansion of Asymptotic Expansion~\ref{prop:nonlinear}, Lemma~\ref{lem:girsanov_bl}, and Proposition~\ref{prop:shifted_domain}, which express both the $O(\sqrt{h})$ and $O(h)$ bias coefficients of general one-dimensional nonlinear Langevin equations explicitly in terms of boundary derivatives of the continuum exit-time and splitting-probability functions and show algebraically that for weak-order-2 integrators driven by Gaussian increments (e.g., stochastic Heun) the expansion collapses to the continuum exit time on a domain widened by $\beta\sigma\sqrt{h}$ plus a universal $+\frac{h}{4}$ variance offset (implying a universal $\tau_{\mathrm{bs}}(h) - \tau(0) = \frac{h}{4} + O(h^{3/2})$ residual for barrier-shifted schemes, so that subtracting $h/4$ yields $O(h^{3/2})$ accuracy at zero extra cost); and (iii)~high-statistics OpenMP Fortran~90 benchmarks and Richardson-extrapolated Fredholm quadrature validating the predicted coefficients, the universal $h/4$ barrier-shift residual, the $O(h^{3/2})$ remainder scaling in the resolved asymptotic regimes, and the $O(h^1)$ convergence restored by barrier shifting and Brownian-bridge exit tests. Practitioners using first-order Euler--Maruyama ($p = 1$) should note that Eq.~\eqref{eq:shifted_domain_nonlinear} leaves an additional interior discretization bias $C_{\mathrm{int}}(x_0)\,h$ after barrier shifting---for instance, $C_{\mathrm{int}} \approx -0.80$ in our Ornstein--Uhlenbeck benchmark, which swamps and reverses the sign of the $+\frac{h}{4}$ offset---so the universal $h/4$ subtraction applies specifically to weak-order-2 schemes driven by Gaussian increments. Extending this boundary-layer matching to multidimensional curved domains---where local mean curvature is expected to enter at order $O(h)$---provides a natural direction for future work.

%======================================================================
\section*{Code Availability}
%======================================================================

The OpenMP-parallelized Fortran~90 simulation driver (\texttt{exit\_time\_benchmarks.f90}), the thread-safe \textsc{pcg64-dxsm} random-number modules (\texttt{precision\_kind\_mod.f90}, \texttt{random\_numbers\_mod.f90}), and the companion quadrature and plotting scripts are publicly archived on Zenodo~\cite{Mannella2026Zenodo}.

\section*{Acknowledgments}

The author thanks the Green Data Center of the University of Pisa for providing the computational resources used in this work.

\section*{Declaration of generative AI and AI-assisted technologies in the manuscript preparation process}

During the preparation of this work, the author used Google Gemini to cross-check intermediate asymptotic and algebraic derivations, assist with auxiliary Python verification and plotting scripts, and improve the clarity and presentation of the text. After using this tool, the author independently verified all analytical results, numerical benchmarks, and code, and takes full responsibility for the content of the published article.

%======================================================================


\begin{thebibliography}{99}
\bibitem{Mannella1999}
R.~Mannella,
\emph{Absorbing boundaries and optimal stopping in a stochastic differential equation},
Phys.\ Lett.\ A \textbf{254}, 257 (1999).

\bibitem{BlytheBray2003}
R.~A.~Blythe and A.~J.~Bray,
\emph{Survival probability of a diffusing particle in the presence of Poisson-distributed mobile traps},
Phys.\ Rev.\ E \textbf{67}, 041101 (2003).

\bibitem{Fritz2025}
L.~Fritz,
\emph{First passage problem: Asymptotic corrections due to discrete sampling},
arXiv:2510.10226 (2025).

\bibitem{Lotov1996}
V.~I.~Lotov,
\emph{On some boundary crossing problems for Gaussian random walks},
Ann.\ Probab.\ \textbf{24}, 2154 (1996).

\bibitem{Chernoff1965}
H.~Chernoff,
\emph{Sequential tests for the mean of a normal distribution IV (discrete case)},
Ann.\ Math.\ Stat.\ \textbf{36}, 55 (1965).

\bibitem{Siegmund1979}
D.~Siegmund,
\emph{Corrected diffusion approximations in certain random walk problems},
Adv.\ Appl.\ Probab.\ \textbf{11}, 701 (1979).

\bibitem{Siegmund1985}
D.~Siegmund,
\emph{Sequential Analysis: Tests and Confidence Intervals}
(Springer, New York, 1985).

\bibitem{BroadieGlassermanKou1997}
M.~Broadie, P.~Glasserman, and S.~Kou,
\emph{A continuity correction for discrete barrier options},
Math.\ Finance \textbf{7}, 325 (1997).

\bibitem{Strittmatter1987}
W.~Strittmatter,
University of Freiburg Report No.\ THEP 87/12, 1987 (unpublished).

\bibitem{Gobet2000}
E.~Gobet,
\emph{Weak approximation of killed diffusion using Euler schemes},
Stoch.\ Process.\ Appl.\ \textbf{87}, 167 (2000).

\bibitem{GobetMenozzi2004}
E.~Gobet and S.~Menozzi,
\emph{Exact approximation rate of killed hypoelliptic diffusions using the discrete Euler scheme},
Stoch.\ Process.\ Appl.\ \textbf{112}, 201 (2004).

\bibitem{GobetMenozzi2010}
E.~Gobet and S.~Menozzi,
\emph{Stopped diffusion processes: boundary corrections and overshoot},
Stoch.\ Process.\ Appl.\ \textbf{120}, 130 (2010).

\bibitem{BouchardGeissGobet2017}
B.~Bouchard, S.~Geiss, and E.~Gobet,
\emph{First time to exit of a continuous It\^o process: General moment estimates and $L_1$-convergence rate for discrete time approximations},
Bernoulli \textbf{23}, 1631 (2017).

\bibitem{JansonsLythe2000}
K.~M.~Jansons and G.~D.~Lythe,
\emph{Efficient numerical solution of stochastic differential equations using exponential timestepping},
J.\ Stat.\ Phys.\ \textbf{100}, 1097 (2000).

\bibitem{Buchmann2005}
F.~M.~Buchmann,
\emph{Simulation of stopped diffusions},
J.\ Comput.\ Phys.\ \textbf{202}, 446 (2005).

\bibitem{Higham2013}
D.~J.~Higham, X.~Mao, M.~Roj, Q.~Song, and G.~Yin,
\emph{Mean exit times and the multilevel Monte Carlo method},
SIAM/ASA J.\ Uncertain.\ Quantif.\ \textbf{1}, 2 (2013).

\bibitem{BernalAcebron2016}
F.~Bernal and J.~A.~Acebr\'on,
\emph{A comparison of higher-order weak numerical schemes for stopped stochastic differential equations},
Commun.\ Comput.\ Phys.\ \textbf{20}, 703 (2016).

\bibitem{Spitzer1956}
F.~Spitzer,
\emph{A combinatorial lemma and its application to probability theory},
Trans.\ Am.\ Math.\ Soc.\ \textbf{82}, 323 (1956).

\bibitem{Feller1971}
W.~Feller,
\emph{An Introduction to Probability Theory and Its Applications}, Vol.~II, 2nd ed.
(Wiley, New York, 1971).

\bibitem{ChangPeres1997}
J.~T.~Chang and Y.~Peres,
\emph{Ladder heights, Gaussian random walks and the Riemann zeta function},
Ann.\ Probab.\ \textbf{25}, 787 (1997).

\bibitem{BallyTalay1996}
V.~Bally and D.~Talay,
\emph{The law of the Euler scheme for stochastic differential equations. I. Convergence rate of the distribution function},
Probab.\ Theory Relat.\ Fields \textbf{104}, 43 (1996).

\bibitem{MilsteinTretyakov2004}
G.~N.~Milstein and M.~V.~Tretyakov,
\emph{Stochastic Numerics for Mathematical Physics}
(Springer, Berlin, 2004).

\bibitem{Mannella2026Zenodo}
R.~Mannella,
\emph{OpenMP Fortran~90 benchmarks for finite time-step bias in stochastic exit problems},
\url{https://doi.org/10.5281/zenodo.24633273}.

\end{thebibliography}
\end{document}